\documentclass[12pt]{amsart}

\usepackage{amsmath,amssymb}

\usepackage{graphicx}

\usepackage{times}

\newcommand{\rmkd}[1]{}

\newcommand{\Z}{{\mathbb Z}}
\newcommand{\Q}{{\mathbb Q}}
\newcommand{\R}{{\mathbb R}}
\newcommand{\g}{\hat{g}}

\newtheorem{thm}{Theorem}
\newtheorem{prop}[thm]{Proposition}
\newtheorem{cor}[thm]{Corollary}

\newcommand{\vett}{v_{*}}

\author{Laura Fainsilber}
\address{Department of Mathematics\\ 
Chalmers University of Technology\\
SE-412 96 Gothenburg  \\
Sweden}

\email{laura@math.chalmers.se}

\author{P\"ar Kurlberg}
\email{kurlberg@math.chalmers.se}

\author{Bernt Wennberg}
\email{wennberg@math.chalmers.se}

\thanks{We acknowledge partial support from the
  National Science Foundation (DMS 0071503) (P.K.), the Swedish Research
  Council  (P.K. \& B.W.), the Royal Swedish
Academy of Sciences (P.K.), and from the EC funded 
RTN network HYKE,  Contract Number : HPRN-CT-2002-00282  (B.W.) }

\title[Lattice points on circles and the Boltzmann equation]%
    {Lattice points on circles and discrete velocity models for the
  Boltzmann equation} 

\begin{document}

\maketitle

\section{Introduction}

The phase space density $f$ of a dilute gas evolves according to the
Boltzmann equation. In the physically relevant case, the gas would be
confined to a subset $\Omega\subset\R^{3}$, and then $f(x,v,t) :
\Omega\times\R^{3}\times\R^{+} \rightarrow \R^{+}$, where  $x$ denotes a
position in space, $v\in\R^{3}$ is a velocity, and $t$ denotes the
time. From a mathematical point of view, it is equally natural to
consider the Boltzmann equation in any spatial dimension, and in some
cases because of 
symmetries of $\Omega$, it is also relevant to consider
$\Omega\subset\R^{d_{1}}$ and $v\in \R^{d_{2}}$ with $d_{1}<d_{2}$.

By a dilute gas we mean one where the particles interact with each
other essentially only by {\em pairwise} interactions. Moreover, the
Boltzmann equation assumes that the particles are so small compared to
other distances, that  they can be considered to be points. 

Under these hypothesis, one can formally derive the Boltzmann
equation~(see~\cite{CIP})
\begin{eqnarray}
  \label{eq1}
  \partial_{t} f(x,v,t) + v\cdot\nabla_{x}f(x,v,t) &=& Q(f,f)(x,v,t)\,.
\end{eqnarray}
The left hand side describes the evolution of the density by free
transport, and the right hand side describes the impact of
collisions. Per definition, a collision is a pairwise interaction that
takes place instantaneously and at one single point in space. Hence
$x$ and $t$ appear only as parameters in $Q(f,f)$, and we
can write
\begin{eqnarray}
\label{eq2}
  Q(f,f)(v) &=& \int_{\R^{d}}\int_{S^{d-1}} \left(
  f(v')f(\vett')-f(v)f(\vett) \right) q(|w|,\cos\theta)
  \,dS(u)  d\vett\,, \nonumber \\
\end{eqnarray}
where the velocities ``before and after a collision'' are related by
\begin{eqnarray}
\label{eq3}
  \begin{array}{lll}
  v' &=&  \frac{1}{2}(v+\vett) + |w|\, u \\ & &\\
  \vett' &=& \frac{1}{2}(v+\vett) - |w|\,  u\,,
  \end{array}
\end{eqnarray}
with $w = (\vett-v)/2$, and with $\cos\theta = \frac{u\cdot w}{|w|}$;
$d\vett$ is the Lebesgue measure in $\R^{d}$, and $dS(u)$ is
the surface measure on $S^{d-1}$. Note
that the  pair of velocities before a collision, $v$ and
$\vett$, and the pair of velocities after the collision,  $v'$ and
$\vett'$, are the endpoints of a diameter on the sphere  which
has its 
centre at $\frac{v+\vett}{2}$ and diameter $|\vett-v|$. This is
exactly the condition needed in order that the collisions preserve the
momentum and energy of the pair of particles. For $d=2$, the sphere
becomes a circle, and this motivates the title of the paper. 

In a discrete velocity model (DVM), the velocities are concentrated on 
a (usually finite)
set of points $v_{j}\in\R^{d}$
in the velocity space: 
\begin{eqnarray*}
  f(x,v,t) &=& \sum_{j} f_{j}(x,t) \delta_{v=v_{j}}\,.
\end{eqnarray*}
The Boltzmann equation~(\ref{eq1}) is then changed into a nonlinear system
of conservation laws,
\begin{eqnarray}
\label{eq4}
  \partial_{t} f_{j} + v_{j}\cdot\nabla_{x}f_{j} =
   \sum_{k,k',j'}  \Gamma_{j,k}^{j',k'} \left( f_{j'}f_{k'} -
   f_{j}f_{k}\right)\, ,
\end{eqnarray}
where the constants $\Gamma_{j,k}^{j',k'}\ge 0$ must be chosen so
that~(\ref{eq4}) makes sense from a physical point of view. In
particular we require that $(v_{j},v_{k})$ and $(v_{j'},v_{k'})$
define two diameters on the same sphere, just as for the usual
Boltzmann equation.  

The first example of a discrete velocity model is that of
Carleman~(\cite{Carleman}), which has two velocities in
$\R$. Many other models have been proposed,  and there is a large
literature on how to construct and analyse 
physically realistic models (i.e., that satisfy the right conservation laws
and an entropy principle), see eg.\cite{BobCer}.

Besides offering many interesting mathematical challenges  (for
example, there is no general theory of global 
existence of solutions to systems like~(\ref{eq4})) the  DVM:s are also
candidates for the numerical approximation of the real Boltzmann
equation~(\ref{eq1}). This leads naturally to the following
question, which is the subject matter of the paper:

Suppose that we choose the discrete set of velocities to be $h
\Z^{d}$, i.e. the integer lattice in $\R^{d}$, scaled by a factor $h$, 
and that we take
\begin{eqnarray*}
  f^{h}(v) = \sum_{\xi \in \Z^{d} } f_{\xi,h} \delta_{v=h\xi}\,,
\end{eqnarray*}
so that $f^{h} \rightarrow f$, in some suitable sense, where $f\in
L^{1}(\R^{d})$. Is it then true that $Q(f^{h},f^{h})(v) \rightarrow
Q(f,f)(v)$ for all $v \in h \Z^{d}$ when $h \rightarrow 0$?

This property, which  is called
{\em consistency}, together with {\em stability} are main ingredients
when proving that a numerical method converges.

The answer is yes. This was proven by Bobylev, Palczewski and
Schneider (\cite{BPS}) for dimensions $d\ge 3$. In this paper,
we prove that it is also true for $d=2$, and hence for all
relevant cases. 

Results of this kind are interesting, because they provide examples
that are relevant to previous results of Desvillettes and
Mischler~(\cite{DesMis}), who proved that solutions to families of
DVM:s can converge  to DiPerna-Lions'  solutions 
to~(\ref{eq1}) if 
certain conditions are satisfied.

Our result  should not, however, be considered as relevant for  numerical
analysis, because the 
rate of convergence is so slow that a numerical method based on the
theory presented here would hardly ever become useful.

The family of models considered here can be seen as coming from a rather
straightforward discretization of the collision integral~(\ref{eq2}).
This integral should be interpreted as an average over the
$2d-1$-dimensional manifold defined by 
\begin{eqnarray}
\label{eq45}
  \lefteqn{\mathcal{M}_{v} = \big\{ (\vett,v',\vett') \in \R^{3 d} \;\;
  \mbox{s.t.} \; \;
  v'+\vett' -\vett = v\, }&&\\ \nonumber
  & & \rule{4cm}{0pt}  |v'|^{2}+|\vett'|^{2} -|\vett|^{2} =
  |v|^{2}   \big\}\,,
\end{eqnarray}
and~(\ref{eq2}) is an iterated integral over this manifold.
For a fixed $v$, we write $w= (\vett-v)/2$, and then~(\ref{eq3}) becomes
\begin{eqnarray*}
  v' &=& v + w +|w|u \\
  \vett' &=& v + w - |w|u\\
\end{eqnarray*}
and also $\vett = v + 2w$. We then write
\begin{eqnarray}
\label{eq45.5}
  g_{v}(w,u) &=& \left( f(v')f(\vett')-f(v)f(\vett) \right)
  q(|w|,\cos\theta)\,, 
\end{eqnarray}
and so (after changing variables in the integral),
\begin{eqnarray*}
  Q(f,f)(v) &=& 2^{d} \int_{\R^{d}}\left(  \int_{S^{d-1}} g_{v}(w,u)\,
  dS(u)\right) \,dw. 
\end{eqnarray*}
If $g$ is sufficiently regular (continuous), and decays sufficiently
rapidly for large $w$, then the Riemann sum for the outer integral converges:
\begin{eqnarray}
\label{eq5}
\lefteqn{ (2h)^d \sum_{\zeta \in \Z^d }
\int_{S^{d-1}}
   g_{v}(h \zeta,u) 
\, dS(u)}\hspace{4em}&&\nonumber  \\
&\longrightarrow & 2^{d}\int_{\R^N} \left( \int_{S^{d-1}}
g_{v}(w,u) \, dS(u) \right) \,dw
\end{eqnarray}
when $h\rightarrow 0$. In order to construct a consistent DVM, it is
then sufficient to evaluate the inner integral in terms of the values
of $g$ on the lattice points $h\Z^{d}$, in such a way  that the result
converges  to $ \int_{S^{d-1}}
g(w,u) \, dS(u) $. While with  the formula~(\ref{eq3}),
the collision integral should be taken over all  
$u\in S^{d-1}$, we have here only access to those $u$ for which
$v'$ and $\vett'$ belong to $h\Z^{d}$. But this is automatically
achieved 
if $\zeta \in  \Z^{d}$,
and if 
$u=\zeta'/|\zeta'|$, where $\zeta'
\in  \Z^{d}$ and $|\zeta'| = |\zeta|$; then for all $v\in h\Z^{d}$,
\begin{eqnarray*}
   v + h \zeta \pm  h|\zeta|u   \in h\Z^{d}\,.
\end{eqnarray*}
However, note that with this construction, the center of the sphere is
restricted to lie on a lattice point,
and so it excludes cases like 
$v=(0,0)$, $\vett=(h,h)$. 

Giving  all points on the sphere  equal weight, one arrives at the
expression 
\begin{eqnarray}
\label{eq6}
 \frac{1}{r_d(|\zeta|^{2})}
\sum_{
\substack{
\zeta' \in \Z^d \\ |\zeta'| = |\zeta|} }
\left( f(v')f(\vett')-f(v)f(\vett) \right) q(|h\zeta|,\cos\theta)\,,
\end{eqnarray}
for approximating the inner integral in~(\ref{eq5}). The function
$r_{d}(n)$ denotes the number of points with integer coordinates on a
sphere in $\R^{d}$ with center at the origin and radius
$\sqrt{n}$, i.e. the number of integer solutions to
$x_{1}^{2}+\cdots+x_{d}^{2}=n$. 

We write, for all $v\in h\Z^{d}$.
\begin{eqnarray}
  \label{eq6a}
  \nonumber
  \lefteqn{ Q^{h}(f,f)(v) =}  && \\
  && (2h)^d \sum_{\zeta \in \Z^d }
  \frac{1}{r_d(|\zeta|^{2})} \sum_{
\substack{
\zeta' \in \Z^d \\ |\zeta'| = |\zeta|} }
\left( f(v')f(\vett')-f(v)f(\vett) \right) q(|h\zeta|,\cos\theta)\,.
\end{eqnarray}

In the two-dimensional case, all the terms in the sum are
$2\pi$-periodic functions of $\theta$, 
and assuming sufficient regularity, they can be expressed as a
convergent Fourier series. It is then  natural to introduce the
 exponential sum
\begin{eqnarray}
\label{eq6b}
S(n,k) = 
\sum_{u \in \Z^2 : |u|^2 = n}
e^{ i k \theta_u}
\end{eqnarray}
where $\theta_u$ is defined by $u =  |u| \cdot  (\sin \theta_u,  \cos
\theta_u)$. 
We will see in Section~\ref{sec:proof}
that to prove that~(\ref{eq6}) converges to the angular integral
in~(\ref{eq5}),
it is enough to prove that for $k\ne0$,
 the terms $ S(n,k)$ converge to zero sufficiently fast as
$n\rightarrow\infty$. Similar exponential sums are relevant for any
dimension, and the work of Bobylev et al. also involves such estimates.

Here  the needed estimate is given as
Proposition~\ref{prop7} in Section~\ref{sec:equid-latt-points}.
Then in Section~\ref{sec:proof} we put the estimates togheter to a
proof of the main result:

\begin{thm}
  \label{th:main:int}
  Consider the Boltzmann equation in two dimensions.
  Assume that $f$ and $q$ are so smooth that the function
  $g_{v}(w,u)$ defined in~(\ref{eq45.5}) is a $C^{2}$-function. Then
  for all $v\in h\Z^{2}$
  \begin{eqnarray*}
    \left| Q(f,f)(v) - Q^{h}(f,f)(v) \right| \rightarrow 0
  \end{eqnarray*}
  when $h\rightarrow 0$.
\end{thm}

Section~\ref{sec:remarks}, finally, contains a numerical illustration on
the distribution terms where the circles contain many points, and we
indicate how the computational cost could be reduced without loosing 
accuracy (still without the intention of actually giving an effective
algorithm).

A more general construction of  discrete velocity
models on scaled integer lattices $h\Z^{2}$ consists in finding sets
of integer points on the manifold $\mathcal{M}$ defined
in~(\ref{eq45}). In this way, mass  and energy conservation are
automatically satisfied, but one also needs to verify that these are
the {\em only} conserved quantities.
 And finally,
in order that the models converge 
to the continuous model when $h\rightarrow 0$, it is necessary that
the integer points are more or less uniformly distributed on
$\mathcal{M}$. 

The models studied here are constructed by
discretizing, one at a time, the iterated integrals~(\ref{eq2}). An
alternative way of writing this integral was introduced by
Carleman~\cite{Carleman}. Using that $v'-v$ and $\vett'-v$ are
orthogonal, one can write (here we specialize to $d=3$)
\begin{eqnarray*}
  \lefteqn{Q(f,f)(v) =} \\
  &&\int_{\R^{3}} \int_{E_{v,v'}}
  \left( f(v')f(\vett')-f(v)f(\vett)  \right) q(w,\cos\theta)
  \frac{1}{|v-v'|^{2}} \,dE(\vett')\,dv'\,,
\end{eqnarray*}
where $E_{v,v'}$ is the plane that contains $v$ and is orthogonal to
$v'-v$, and where $dE(\vett')$ is the Euclidean measure on this
plane. Heintz and Panferov~\cite{HeinPanf} have analysed a DVM based
on this interated integral, and proved that the method is consistent
with the continuous model. This is somehow easier, because on {\em
  all} planes, the integer points are uniformly distributed, and they
are all found by solving linear Diophantine equations. However, the
density of points depends strongly on $v'-v$, and so it is far from
trival to prove the consistency. And again, the two-dimensional
situation is more difficult, and has not yet been studied.

Yet another approach was introduced by Rogier and Schneider~\cite{RS},
who used the
theory of Farey series to discretize the angular variable in the
collision integral.

\section{Number theoretic background}
\subsection{Points on spheres; Asymptotics}

To prove that~(\ref{eq6}) converges to the correct limit when
$h\rightarrow 0$, one has to study the set 
\begin{eqnarray*}
\{ \zeta/|\zeta| : \zeta \in \Z^d, |\zeta|^2 = n \}
\end{eqnarray*}
and to show that the points of this set are sufficiently well
distributed on $S^{d-1}$ when $n$ is large; it is here that the number
theoretical issues enter the game.
\newcommand{\C}{{\mathbb C}}
Indeed, we can view the set of points with integer coordinates on a sphere 
of squared radius $n$ centered at the origin,
$$ 
\{ (x_1, \dots x_d) \in \Z^d, \sum_{i=1}^d x_i^2=n\}\,,
$$
 as the  solution set for a quadratic form, 
and use the theory of integral quadratic forms to get 
estimates on the number of points (see for instance \cite{GRO}).
The expected number of points with integer coordinates on a
sphere 
clearly depends
on the dimension $d$. 
The naive approach to find 
the order of magnitude
for a given dimension is 
to use the volume of a ball, divided by the number of spheres contained 
in the ball.  The volume of a ball of radius $\sqrt n$ grows as 
$n^{d/2}$ whilst the number of spheres is $n$.
For $d=2$, this leads 
us to expect a constant number of lattice points on circles, for $d=3$ a  
growth proportional to $\sqrt n$, etc.
However, for small $d$ this approach is misleading; 
the growth is quite irregular, and depends on the divisor structure of
$n$.   
For $d=2$, we will see below that only values of $n$ of the form
$n=2^s q^2 p_1^{\alpha_1}\dots p_r^{\alpha_r}$,
where $q$ is a product of primes of the form $4k+3$ and
the $p_i$'s are primes of the form $4k+1$ (see below),
yield circles with lattice points, and thus most circles 
have no points at all.  In fact, Landau proved in 1908 that the number
of circles with at least one lattice point,  
of integer squared radius 
smaller than $x$, grows as $C x/\sqrt{\log x}$.
Moreover, there are also infinite families of circles
with very few lattice points; radii that are a power of 2 yield 4 
points for instance, and
radii that are the square root of a prime of the form $p=4k+1$ yield exactly
8 points. On the other hand,
the number of lattice  points on a circle is not bounded, for instance
a circle with $n= p_1\dots p_r$ as above where all the $p_i$ are distinct from 
each other has $4\cdot 2^r$ points.

In dimension 3, all values of $n$ not of the form $n=4^s(8k+7)$ yield 
spheres containing points with integer coordinates.  
This still leave a fairly large number of spheres with no points on,
but for our purposes this does not really matter, as 
such spheres 
do not appear in the summation formulas (there is no relevant value 
for $\zeta$.)  Among the spheres with lattice points, multiplying the
radius by a  
power of 2 does not increase the number of points, but if we correct 
for this fact, the ratio between the number of points and 
the naive estimate is bounded, up to constants only depending on
$\epsilon$, from above by  
$n^{ \epsilon}$, and below by $n^{-\epsilon}$
for all $\epsilon>0$ (see \cite{GRO} 
Ch.~4 for exact formulas involving class numbers or $L$-series.)

The higher-dimensional cases behave in a somewhat more regular fashion.
Lagrange proved that every
positive integer can be written as the sum of four squares, and thus for 
dimension $d \geq 4$, every sphere whose squared
radius is an integer, has lattice points. 
For $d=4$ the number of points still oscillates rather wildly, with 
spheres with radius a power of $2$ having just $24$ points, but for 
greater dimensions, the naive estimate gives  the
correct asymptotic growth of the number of points.

Getting circles (or spheres) with ``sufficiently many'' lattice
points, however, is not quite enough for our purposes: we also need
that the lattice points are sufficiently uniformly distributed when
projected on the unit sphere.
In dimensions $3$ and higher, this 
follows from estimates on Fourier coefficients of modular
forms.  The case $d\geq 4$, with some restrictions on the set of
numbers in which $n$ tends to infinity when $d=4$, is due to
Pommerenke~\cite{pommerenke}.  For $d=3$,
Duke~\cite{duke} and Golubeva-Fomenko~\cite{golubeva-fomenko}
used Iwaniec's~\cite{iwaniec} estimates
on Fourier coefficients of half integral weight forms
to obtain uniform distribution.
Unfortunately, these techniques do not apply in dimension $2$.  
Moreover, 
there are circles with large number of lattice points that are poorly 
distributed:  
\begin{thm}{\rm (Cilleruelo \cite{CIL})}
For any $\epsilon >0$ and for any integer $k$, there exists a circle
$x^2+y^2=n$  
with more than $k$ lattice points such that all the lattice points are
on the arcs 
$\sqrt n e^{(\pi/2)(t+\theta)i}$ with $|\theta|< \epsilon$, $t\in\{0,1,2,3\}$.
\end{thm}
On the other hand, we may use some other techniques
from analytic number theory
to show that lattice points on circles are equidistributed {\em on
average}, and this is good enough for our purpose.

\subsection{From points on circles to Gaussian integers}
\label{sec:gaussian-integers}

In the plane, we can view lattice
points on a circle of  radius $\sqrt{n}$, centered at the origin, 
as complex numbers with integer  real and imaginary parts,
and squared modulus~$n$.  
It might seem as a trivial restatement, but doing so allows us to use
use some  techniques from algebraic number theory.
The {\em Gaussian integers}, i.e., the set
$$\Z[i]=\{x+iy \in \C, (x,y)\in \Z^2\},$$
is the ring of integers of the field $\Q(i)$.  It shares an
important 
property with the ordinary integers, namely unique
factorization\footnote{This is rather unusual, the ring of integers in
  most number fields will not have   this property.},
i.e., just as every integer in $\Z$ factors into prime numbers, and the 
factorization is unique  up
to ordering the primes and multiplying by $-1$, Gaussian integers factor 
into Gaussian primes, uniquely up to ordering and multiplication by 
$-1, i, -i$ (these and $1$ are the units, i.e. the elements having a 
multiplicative inverse in $\Z[i]$). 
For a more thorough introduction to primes in quadratic number fields, 
see for instance~\cite{hardy-wright-third}, Ch. XV.

The Gaussian primes (i.e. the elements of $\Z[i]$ that cannot be written as 
a product of Gaussian integers with smaller modulus), are of three types: 
\begin{itemize}
 \item the prime numbers $q\in\Z$ such that $q\equiv 3 \mod 4$ remain
prime in $\Z[i]$ (e.g. 3, 7, 11, 19,...);

\item for prime numbers $p\in\Z$ such that $p \equiv 1 \mod 4$, there 
exists $x,y \in \Z$ s.t. $p= x^2+y^2$.  Hence $p$ factors in 
$\Z[i]$ as a product of two Gaussian primes
$$
p = (x+iy)(x-iy)
$$
(e.g. $5$ factors into $(2+i) (2-i)$ in $\Z[i]$)
\item last (and least!), $1+i$ is prime
(note that  $(1+i)(1-i)=2$ and that
$1-i=-i(1+i)$ is merely ``another form of the same prime''
just as $3$ and $-3$ represent the same prime).
\end{itemize}

\bigskip

If $n$ is the sum of two squares, then it can be factored in
$\Z[i]$:
\begin{eqnarray*}
  n = X^{2} + Y^{2} = ( X + i Y ) (X - iY).
\end{eqnarray*}
If $z=x+iy$ is a prime factor
of $ X + i Y $, then $\bar{z}=x-iy$ must be a prime factor of $ X - i
Y $. It follows that prime factors $q\equiv 3 \mod 4$ of $n$ must
appear in even powers. In addition, multiplying $n$ by an even power
of a prime $q$ that is congruent with $ 3 \mod 4$ changes neither the
number of solutions to  $n = X^{2} + Y^{2}$ nor the distribution of
arguments of the solutions. 

Suppose now that $n$ contains a factor $p^\alpha $, where $p\equiv 1 \mod
4$. The number $p$ can be factored in $\Z[i]$ as $(x+iy)(x-iy)$, and
hence the 
multiplicity of $x+iy$ as a factor of $n$ is $\alpha$, and the same is
true for $x-iy$. It follows that the multiplicity of $x+iy$ in $X + i
Y$ can be any integer $j$, with $0\le j\le \alpha$, and the 
multiplicty of $x-iy$ is  then $\alpha-j$.

The same calculation can be done for powers of 2; however, the solutions
given by different choices of $j$ in that case differ by a
multiplication by a power of $i$, and so the power of 2 
does not influence the number of solutions.

All  solutions to $n= X^{2} + Y^{2}$ can now be expressed as
$X+iY = \sqrt{n} \exp(i \theta)$, where all possible values of the
argument $\theta$ can be computed as sums of terms deriving from the
different factors of $n$ in the following way:

\begin{enumerate}
\item 
$X+iY$ can be multiplied by any unit, i.e. by $\pm 1$ or $\pm
  i$. This gives a term $k \pi/2$ in the argument, $k=0,1,2,3$. 
\item If the multiplicity of 2 in $n$ is odd, then the argument must contain
  $\pi/4$, the argument of   $1+i$; the number of solutions does not
  change.
\item For each prime factor $p\equiv 1 \mod 4$
  in $n$, let  $\alpha_p$ be the multiplicity of $p$ in $n$,  let $p = x_p^2 +
  y_p^2$, and set $\theta_p =   \arg(x_p+iy_p)$. For a particular choice of
  $j$, $0\le j\le \alpha_p$, the   argument added to $X+iY$ is $j \theta_p -
  (\alpha_p-j)\theta_p = (2j-\alpha_p)   \theta_p$.  
\end{enumerate}
Since the choices of $k$, and of the different $j's$ are independent,
the number of different solutions is $\quad 4 \displaystyle
 \prod_{p\equiv1\!\!\!\mod\!4}
(\alpha_p+1)$.

\bigskip

\subsection{Results on the distribution of primes and on the angular distribution of points}
\label{sec:some-results-from}
We will need the following results:
 
\begin{thm}[Merten's Theorem, see \cite{hardy-wright-third}, Ch. 22.8] 
$$
\prod_{\substack{p \leq x \\ \text{$p$ prime}}} (1-1/p) \sim
e^{-\gamma}/\log x\,,  
$$
where $\gamma \simeq 0.57$ is Euler's constant.
\end{thm}

As for the angular distribution of Gaussian primes, a result by
  Kubilyus gives that the angles $\{\theta_p\}_{p
  \equiv 1 \mod 4}$ are equidistributed in $[0,\pi/4]$ in the
  following sense:

\begin{thm}[Kubilyus, \cite{kubilyus-gaussian-primes-in-sectors}]
\label{thm:kubi}
The number of Gaussian primes $\omega$ in the sector $0 \leq \alpha
\leq \arg(\omega) \leq \beta \leq 2 \pi$, $|\omega|^2 \leq u$ is equal
to 
$$
\frac{2}{\pi}(\beta-\alpha) \int_{2}^u \frac{dv}{\log v} + 
O\left( u \exp(-b \sqrt{\log u}) \right)\,,
$$
where $b$ is an absolute positive constant.
\end{thm}
From Kubilyus' Theorem, it is straightforward to deduce (see
\cite{erdos-hall-gaussian-integers}, p.~92):
\begin{cor}
\label{cor:gaussian-merten}
If  $k \in 4 \mathbb{N}$ and $\log k \leq b \sqrt{\log x}$, then
$$
\sum_{\substack{p \leq x\\p \equiv 1 \mod 4}}
\frac{|\cos(k \theta_p)|}{p}
\leq 
\frac{1}{\pi} \log \log x +
(1-2/\pi) \log \log k + O(1).
$$  
\end{cor}

\section{Equidistribution of lattice points on circles}
\label{sec:equid-latt-points}

What is needed for the proof of consistence of the discrete velocity
model are estimates on the equidistribution of lattice points on circles.

The aim of this section is to show that  lattice
points on circles are equidistributed {\em on average} in the
sense that the exponential sums $S(m,k)$ converge to zero when $m$
goes to infinity.
We recall the definition of $S(m,k)$:
\begin{eqnarray*}
  S(m,k) = \sum_{|w'|^{2} = m} e^{ik \theta_{w'}}\,.
\end{eqnarray*}

\begin{prop}
\label{prop7}
If $4 \nmid k$ then $|S(m,k)| = 0$.  
If $4|k$ and $k \neq 0$, there exist $C$ and $b>0$ such that 
$$
\log\left( \frac{1}{X}\sum_{m \leq X}
|S(m,k)|\right)
\leq
C -   (1- 2/\pi) \log \left( \frac{\log X}{\left(\log |k|\right)^{2}}  \right)  
$$
for  $X$ sufficiently large and 
$\log|k| \leq b \sqrt{\log X} $.
\end{prop}

{\em Remark:} 
The mean discrepancy of the angles of Gaussian integers were studied
by K{\'a}tai and K{\"o}rnyei in
\cite{katai-kornyei-lattice-points-on-circles}, and by 
Erd{\H{o}}s and Hall in
\cite{erdos-hall-gaussian-integers}.  Our method is similar
to theirs, except that they  bound 
$$
\frac{1}{X/\sqrt{\log X}}\sum_{m \leq X}
\frac{|S(m,k)|}{r(m)}
$$
instead of 
$$\frac{1}{X}\sum_{m \leq X} |S(m,k)|.$$

The proof is based on the observation that $|S(m,k)|/4$ is a
{\em multiplicative} function, i.e. a function 
 $f:\mathbb{N} \to \C$ such that 
$f(mn)=f(m)f(n)$ for all $m,n$ such that $(m,n)=1$.
It turns out that the mean value of a multiplicative function, under
fairly general circumstances, can be bounded in terms of 
an exponential of a sum  over primes.  
To make the paper more self contained, we include a weak form of the 
{\em Halberstam-Richert inequality}
(cf. \cite{halberstam-richert-inequality}).

\begin{thm}
\label{thm:multiplicative-average}
Let $f$ be a nonnegative multiplicative function such that 
\begin{equation}
  \label{eq:4:9}
\sum_{n \leq
  x} f(n) = O(x)\,,
\end{equation}
and $f(p^k) = O(k)$ for all primes $p$ and $k \geq 1$.
Then there exists $C>0$ such that
$$
\frac{1}{X}
\sum_{m \leq X}
f(m)
\leq C \cdot
\exp \left( \sum_{p\leq X}  \frac{f(p)-1}{p} \right) 
+ O(\frac{1}{\log X})
$$
for all sufficiently large $X$.
\end{thm}
\begin{proof}
Following Wirsing \cite{Wirsing}, let 
$$
F(t) = \sum_{n \leq t} f(n)\,.
$$
Then 
$$
\int_1^X \frac{F(t)}{t} \, dt
= 
F(X) \log X + O(1) -
\sum_{n \leq X} f(n) \log n\,.
$$
On the other hand, by assumption, we have $F(t) =O(t)$,
thus 
$$
\int_1^X \frac{F(t)}{t} \, dt
=O(X)\,,
$$
and hence
$$
F(X) \log X 
\leq
O(1) + X +
\sum_{n \leq X} f(n) \log n\,.
$$
Using $\log n = \sum_{d|n} \Lambda(d)$, where $\Lambda$ is the
von Mangoldt function\footnote{That is,  $\Lambda(d) = \log p$ if
  $d=p^k$ and $k \geq 1$, otherwise $\Lambda(d) = 0$.}
we have 
$$
\sum_{n \leq X} f(n) \log n
=
\sum_{n \leq X} f(n)
\sum_{d|n} \Lambda(d)
=
\sum_{d \leq X} 
\Lambda(d)
\sum_{m \leq X/d} 
f(dm)
$$
\begin{equation}
  \label{eq:4:1}
=
\sum_{d \leq X} 
\Lambda(d)
\sum_{\substack{m \leq X/d, \\ (m,d)=1}}
f(dm)
+
\sum_{d \leq X} 
\Lambda(d)
\sum_{\substack{m \leq X/d, \\ (m,d)>1}}
f(dm)\,.
\end{equation}

Now, since $\Lambda(d) = 0$ unless $d$ is a prime power, we 
have 
\begin{multline}
\label{eq:4:2}
\sum_{d \leq X} 
\Lambda(d)
\sum_{\substack{m \leq X/d, \\ (m,d)>1}}
f(dm)
=
\sum_{\substack{p^{k+l} \leq X\\k,l \geq 1}} 
\log(p)
\sum_{\substack{m \leq X/p^{k+l} \\ (p,m)=1}}
f(p^{k+l}m)
\\
=
\sum_{\substack{p^{k+l} \leq X\\k,l\geq 1}} 
\log(p) f(p^{k+l})
\sum_{\substack{m \leq X/p^{k+l} \\ (p,m)=1}}
f(m)\,.
\end{multline}
By the assumptions on $f$, 
$$
f(p^{k+l})\sum_{\substack{m \leq X/p^{k+l} \\ (p,m)=1}} f(m)
\leq
O(k+l)
\sum_{\substack{m \leq X/p^{k+l} }} f(m)
=O \left( (k+l) \frac{X}{p^{k+l}} \right)\,,
$$
and thus the second term in (\ref{eq:4:1}) is 
$$
=
O \left(
\sum_{\substack{p^n \leq X\\n \geq 2}} 
\log(p) n^2
\frac{X}{p^{n}}
\right)
=O(X)\,,
$$
since 
$$
\sum_p
\sum_{n\geq 2} 
\log(p) n^2 p^{-n} 
\leq 
\sum_p
\frac{ \log(p)}{p^{2}} \sum_{m \geq 0} (2+m)^2 2^{-m}
< \infty\,.
$$

As for the first term in (\ref{eq:4:1}), we have (recall that $f$ is
multiplicative and nonnegative) 
$$
\sum_{d \leq X} 
\Lambda(d)
\sum_{\substack{m \leq X/d, \\ (m,d)=1}}
f(dm)
=
\sum_{d \leq X} 
\Lambda(d) f(d)
\sum_{\substack{m \leq X/d, \\ (m,d)=1}}
f(m)
$$

$$
\leq
\sum_{m \leq X} 
f(m)
\sum_{d \leq X/m} 
\Lambda(d)
f(d)\,.
$$
Now,
$$
\sum_{d \leq X/m} 
\Lambda(d)
f(d)
=
\sum_{\substack{p^k \leq X/m \\ k \geq 1}} 
\log(p)
f(p^k)
\leq 
\sum_{\substack{p^k \leq X/m \\ k \geq 1}} 
\log(p)
O(k)
=
O(X/m)
$$
since
$$
\sum_{p \leq X/m}
\log(p)
= O(X/m)
$$
by the Prime number theorem, and 
$$
\sum_{\substack{p^k \leq X/m \\ k \geq 2}} 
k\log(p)
= O \left( (X/m)^{1/2} \log^3(X/m) \right)
= O( X/m)\,.
$$
Thus,
$$
\sum_{m \leq X} 
f(m)
\sum_{d \leq X/m} 
\Lambda(d)
f(d)
= O\left(
\sum_{m \leq X} 
f(m) \frac{X}{m}
\right)\,.
$$
But since $f$ is nonnegative and multiplicative, we have
$$
\sum_{m \leq X} 
 \frac{f(m)}{m}
\leq
\prod_{p \leq X}
\left( 1+f(p)/p + f(p^2)/p^2 + \ldots \right) 
$$
$$\leq
\prod_{p \leq X}
\left( 
\left(  1+f(p)/p \right) \cdot 
\left(1+ f(p^2)/p^2 + f(p^3)/p^3 + \ldots \right)
\right) \,,
$$
and since 
$$
\sum_{p \leq X}
\left(
f(p^2)/p^2 + f(p^3)/p^3 +\ldots \right)
\leq
\sum_{p}
\sum_{k \geq 2 }
\frac{O(k)}{p^k}
< \infty\,,
$$
we find that 
$$
\sum_{m \leq X} 
\frac{f(m)}{m}
= O \left(
\prod_{p \leq X}
\left( 1+f(p)/p \right)
\right)\,.
$$
Thus,
$$
F(X) \log X 
= O \left(
X + X \cdot 
\prod_{p \leq X} \left( 1+f(p)/p \right)
\right)\,,
$$
hence
$$
\frac{F(X)}{X}
= O \left( \frac{1}{\log X}
+ \frac{\prod_{p \leq X} \left( 1+f(p)/p \right)}{\log X}
\right)\,.
$$
Now, by Merten's theorem, we have
$$
\prod_{p \leq X} (1-1/p)
\sim \frac{e^{-\gamma}}{ \log X}\,,
$$
and thus
$$
\frac{F(X)}{X}
=
O \left(
\frac{1}{\log X} +
\prod_{p \leq X} \left( 1+\frac{f(p)-1}{p} -\frac{f(p)}{p^2} \right)
\right)
$$
$$=
O \left(
\frac{1}{\log X} +
\exp \left( 
\sum_{p \leq X}
\frac{f(p)-1}{p} 
\right)
\right)\,.
$$

\end{proof}

\begin{proof}[Proof of Proposition~\ref{prop7}]

To see that $|S(m,k)/4|$ is a multiplicative function, it is enough to
recall the factorization of $m$ into Gaussian primes. Namely, if
$p_{1}^{\alpha_{1}}, ... p_{J}^{\alpha_{J}}$ are all prime factors of $m$ with
  $ p \equiv 1 \mod 4$, 
\begin{eqnarray*}
  S(m,k) = \sum_{\ell=0}^{3}i^{k\ell} 
   \sum_{j_{1} = 1}^{\alpha_{1} } \cdots
   \sum_{j_{J} = 1}^{\alpha_{J}}  e^{ik (\theta_{0} +
   (\alpha_{1}-2j_{1})\theta_{p_{1}} + ... +
   (\alpha_{J}-2j_{J})\theta_{p_{J}}) } \,. 
\end{eqnarray*}
Here $\theta_{0}$ is a multiple of $\pi/4$ which comes from powers of
$2$ in $m$, and the $\theta_{p_{j}}$ can be computed from the Gaussian
factorization as described in Section~\ref{sec:gaussian-integers}. 
Also, because $\sum_{\ell=0}^{3}i^{k\ell} = 4$ if $4 \,|\, k $ and
zero otherwise, 
\begin{eqnarray*}
\frac{|S(m,k)|}{4} = \left|
\sum_{j_{1} = 1}^{\alpha_{1} } \cdots
   \sum_{j_{J} = 1}^{\alpha_{J}}  e^{ik ( 
   (\alpha_{1}-2j_{1})\theta_{p_{1}} + ... + (\alpha_{J}-2j_{J})\theta_{p_{J}}) }\,
\right| \,, 
\end{eqnarray*}
and this sum clearly factorizes, each factor containing a sum of terms
corresponding to one of the prime factors  $p$. Hence
$$
f_k(m) = \frac{|S(m,k)|}{4}
$$
is a nonnegative multiplicative function, as stated. In addition it satisfies
$f_k(m) \leq r(m)/4$ for all $m$.  Thus, since 
$$
\sum_{n \leq T} r(n) = |\{x,y \in \Z : x^2 + y^2 \leq T  \}| 
\sim 
\pi \left( \sqrt{T} \right)^2 = \pi T
$$
we have
$$
\sum_{n \leq T} f_k(n) = O(T).
$$
Moreover, if $p \equiv 3 \mod 4$ then 
\begin{equation}
  \label{eq:4:12}
f_k(p^l) = 
\begin{cases} 
1 \text{ if $l$ is even,}\\
0 \text{ if $l$ is odd,}
\end{cases}
\end{equation}
and if $p \equiv 1 \mod 4$ then 
\begin{equation}
  \label{eq:4:13}
f_k(p^{l}) = \left| \sum_{j=0}^l  e^{i k (l-2j)\theta_p  } \right|\,,
\end{equation}
and thus  $f_k(p^l) \leq l+1$ for all prime $p$ and $l\geq 1$.  The
assumptions in  
in Theorem~\ref{thm:multiplicative-average} are thus satisfied, and we
obtain 
$$
\frac{1}{X}
\sum_{m \leq X}
|S(m,k)| = 
\frac{4}{X}
\sum_{m \leq X}
f_k(m)
\leq C
\exp\left(
\sum_{p\leq X}  \frac{f_k(p)-1}{p} 
\right)
+O \left( \frac{1}{\log X} \right)\,.
$$
Now, by (\ref{eq:4:12}) and (\ref{eq:4:13}), we have 
$$
f_k(p) = 
\begin{cases}
2 | \cos( k\theta_p)| & \text{if $p \equiv 1 \mod 4$,}\\
0 & \text{if $p \equiv 3 \mod 4$}\,.
\end{cases}
$$
Hence 
$$
\sum_{p\leq X}  \frac{f_k(p)-1}{p}
=
\sum_{\substack{p\leq X \\ p \equiv 1 \mod 4}}  
\frac{2 |\cos(k \theta_p)|}{p}
-
\sum_{p\leq X}  
\frac{1}{p}\,.
$$
By Corollary~\ref{cor:gaussian-merten},
$$
\sum_{\substack{p\leq X \\ p \equiv 1 \mod 4}}  
\frac{2|\cos(k \theta_p)|}{p} 
\leq
\frac{2}{\pi} \log \log X +
2(1-2/\pi) \log \log k + O(1).
$$
if $\log k \leq b \sqrt{\log X}$.  
By Merten's theorem, 
$$
\sum_{p\leq X}
\frac{1}{p} = \log \log X + O(1)\,,
$$
and thus 
$$
\sum_{p\leq X}
\frac{f_k(p)-1}{p} 
\leq
(2/\pi-1) \log \log x +
2(1-2/\pi) \log \log k + O(1)\,.
$$
\end{proof}

\section{Proof of Theorem~\ref{th:main:int}}
\label{sec:proof}

Here we carry out the steps of the proof as indicated in the
introduction. First recall that the collision operator can be written 
\begin{eqnarray}
\label{eq:pro:1}
  Q(f,f)(v) &=& 4 \int_{\R^{2}}\left(  \int_{-\pi}^{\pi}  g_{v}(w,\theta)\,
  d\theta\right) \,dw,   
\end{eqnarray}
where, if we identify $u\in S^{1}$ with $\theta\in [-\pi ,\pi[$, 
\begin{eqnarray*}
  g_{v}(w,\theta) = q(|w|,\cos(\theta))\left(
  f(v')f(\vett')-f(v)f(\vett)\right),
\end{eqnarray*}
and 
\begin{eqnarray*}
  v' &=& v + w  +  R_{\theta} w \\
  \vett' &=& v + w  -  R_{\theta} w\,;
\end{eqnarray*}
as before, $w=(\vett-v)/2$, and $R_{\theta}$ denotes a rotation by an
angle~$\theta$. Writing the Boltzmann equation for two-dimensional
velocities, of course we have already stepped away from the physically
realistic case, but disregarding this, a  common  assumption on 
$q$ is that 
\begin{eqnarray*}
   q(|w|,\cos(\theta)) &=& q_{1}(|w|) q_{2}(\theta)\,,   
\end{eqnarray*}
where $q_{1}(|w|) \sim |w|^{\alpha}$ for some $\alpha\in[0,1]$, and
where $q_{2}(\theta) \sim |\theta|^{-\gamma}$ for some
$\gamma\in]1,3[$. This corresponds to a molecular interaction by
hard inverse power law forces. With the stronger assumption that
$q_{1}$ is smooth and strictly positive, it is possible to prove that
there is a smooth solution $f(v,t)$ to the Boltzmann equation~(see
\cite{DeWe}), and then this also gives some regularity to
$g(w,\theta)$, in spite of the singularity of~$q_{2}$.

However, much work on the Boltzmann equation has been done with the
hypothesis that $q$ is bounded or continuous with respect to
$\theta$. With that assumption, the solution $f(v,t)$ keeps exactly
the regularity of the initial data.

Because of this, it is relevant to assume whatever regularity of the
solutions that is needed for the computations. With the aim of 
making the calculations easy, Theorem~\ref{th:main:int} has been
written with  unnecessarily strong hypothesis.

To simplify notation a little, let
\begin{eqnarray*}
  G_{v}(w) &=& \int_{-\pi}^{\pi}
  g_{v}(w,\theta)\,d\theta\,,
\end{eqnarray*}
in the continuous case, and for the discrete case (then we assume, of
course, that $v\in h\Z^{2}$)
\begin{eqnarray*}
   G^{h}_{v}(h\zeta) &=&\frac{1}{r(|\zeta|^{2})} \sum_{\substack{
\zeta' \in \Z^d \\ |\zeta'| = |\zeta|}}
  g_{v}(h\zeta,\theta)\,,
\end{eqnarray*}
where  $\theta$ is the angle between $\zeta'$ and $\zeta$. As before, 
$r(|\zeta|^{2})$ denotes the number of integer points on a sphere with
radius $|\zeta|$.

Let
\begin{equation}
  \label{eq:zhr}
  Z_{h,R} = \{ z\in \Z^{2} \;\;\mbox{s.t.}\;  \; |z| \le R/h \}
\end{equation}
 for some $R>0$
(this is the most natural example, but other choices might be more
efficient, as we shall see later). 
 We want to prove that
\begin{eqnarray}
\label{eq:pr0}
  Q(f,f)(v) - (2h)^{2} \sum_{\zeta \in Z_{h,R}} 
                G^{h}_{v}(h\zeta)  \rightarrow 0
\end{eqnarray}
when $h\rightarrow 0$, and also make as precise a statement as
possible about the rate of convergence.

\begin{thm} 
\label{thm:main}
Suppose that $g_{v}(w,\theta)$ in~(\ref{eq:pro:1})
satisfies
\begin{enumerate}
\item $g_{v}(w,\theta)$ is a $C^{1}$-function w.r.t. $w$
\item $g_{v}(w,\theta)$ is a $C^{2}$-function w.r.t. $\theta$
\item $ \| g_{v}(\cdot,\theta) (1+|\cdot|^{2}) \|_{L^{1}(dw)} \le C $ 
\end{enumerate}
(This holds e.g. if the function $f$ and the crossection $q$ are
$C^{2}$.) For given $R>0$ and $h>0$, let $Z_{h,R}$ be as in~(\ref{eq:zhr}).
Then given $\varepsilon>0$ there are reals $R>0$ and $h>0$ such that
  \begin{eqnarray*}
 \left| Q(f,f)(v) - (2h)^{2} \sum_{\zeta \in Z_{h,R}} 
                G^{h}_{v}(h\zeta) \right|  \le \varepsilon \, .
  \end{eqnarray*}
\end{thm}

\begin{proof}
 We still consider $Q(f,f)$ as an iterated
integral, and write (for $v\in h\Z^{2}$)
\begin{eqnarray}
\label{eq:pr1}
\lefteqn{  Q(f,f)(v) - (2h)^{2} \sum_{\zeta\in Z_{h,R}} 
                G^{h}_{v}(h\zeta)   \nonumber } \hspace{4em} &&  \\
    &=&  \int_{\R^{2}} G_{v}(w) \,dw - 
        (2h)^{2} \sum_{\zeta\in Z_{h,R}} G_{v}(h\zeta)
    \nonumber \\
    &&+  (2h)^{2} \sum_{\zeta\in Z_{h,R}}
         \left(  G_{v}(h\zeta) -   G^{h}_{v}(h\zeta)   \right)\,. 
\end{eqnarray}

From the third part of the hypothesis on $g$  (which is implied by a decay
of $f(v)$ for large velocities),  it follows that for all $R>0$,
\begin{eqnarray}
\label{eq:pr1.5}
  \int_{|w| \ge  R } G_{v}(w) \,dw \le  \frac{C_{1}}{R^{2}}\,.
\end{eqnarray}
Continuity of $G_{v}(w)$ would be enough to conclude  that
\begin{eqnarray*}
  \left| \int_{|w| <  R } G_{v}(w) \,dw - (2h)^{2} \sum_{\zeta\in
  Z_{h,R}} G_{v}(h\zeta)\right|
\rightarrow 0
\end{eqnarray*}
when $h\rightarrow0$. The hypothesis on $g_{v}(w,\theta)$ implies that
actually $G_{v}(\cdot)\in C^{1}$, and there is a constant $C_{2}$ such
that the difference is smaller than 
\begin{eqnarray}
\label{eq:pr2}
 C_{2}R^{2} h = C  \max_{w,j} |\partial_{w_j} G_{v}(w)| \; R^{2} h  \,.
\end{eqnarray}

Next we turn to the difference $ G_{v}(h\zeta) -   G^{h}_{v}(h\zeta)$,
i.e. of
\begin{eqnarray}
\label{eq:pr3}
 \frac{1}{2\pi} \int_{-\pi}^{\pi}
  g_{v}(h\zeta,\theta)\,d\theta -  \frac{1}{r(|\zeta|^{2})} \sum_{\substack{
\zeta' \in \Z^2 \\ |\zeta'| = |\zeta|}}
  g_{v}(h\zeta,\theta)\,,
\end{eqnarray}
(recall that in the second term, $\theta$ is the angle between
$\zeta'$ and $\zeta$).
We first write the periodic function $g_{v}(h\zeta,\theta)$ as a Fourier
series,
\begin{eqnarray*}
 g_{v}(h\zeta,\theta) = 
 \sum_{k \in \Z} \g_{v}(\zeta,k)  e^{i k \theta}\,,
\end{eqnarray*}
where
\begin{eqnarray*}
  \g_{v}(\zeta,k) = \frac{1}{2\pi} \int_{-\pi}^{\pi}
  g_{v}(h\zeta,\theta) e^{-i k 
  \theta}\,d\theta\,.
\end{eqnarray*}
The assumptions on $g$  imply the existence of a constant
$C_{3}$ so that
\begin{eqnarray}
\label{eq:pr3.5}
  | \g_{v}(\zeta,k) | \le \frac{C_{3}}{1+k^{2}} \,.
\end{eqnarray}
Then~(\ref{eq:pr3}) becomes
\begin{eqnarray*}
 \g_{v}(\zeta,0) - \frac{1}{r(|\zeta|^{2})} \sum_{\substack{
  \zeta' \in \Z^2 \\ |\zeta'| = |\zeta|}}
    \g_{v}(\zeta,0)
+ \frac{1}{r(|\zeta|^{2})} \sum_{\substack{
 \zeta' \in \Z^2 \\ |\zeta'| = |\zeta|}} \sum_{k\ne0} 
 \g_{v}(\zeta,k)  e^{i k \theta}\,,
\end{eqnarray*}
where the first terms cancel out, and only last sum remains. We  next
split that sum into a part with 
$|k|\le M$, and a remainder,  which can be made small by choosing $M$
large, if $g$ is sufficiently smooth with respect to
$\theta$. Using~(\ref{eq:pr3.5}), 
\begin{eqnarray*}
 \Bigg| \frac{1}{r(|\zeta|^{2})} \sum_{\substack{
 \zeta' \in \Z^2 \\ |\zeta'| = |\zeta|}} \sum_{|k| \ge  M} 
 \g_{v}(\zeta,k)  e^{i k \theta} \Bigg| \le
    2\frac{C_{3}}{M}\,.
\end{eqnarray*}
To find the contribution of this term to~(\ref{eq:pr1}), we multiply
by $(2h)^{2}$ and sum over $\zeta\in Z_{h,R} $ to find a bound of the
form
\begin{eqnarray}
\label{eq:pr:21}  
 \frac{R^{2}  C_{4}}{M}\,.
\end{eqnarray}

For the remaining part, using~(\ref{eq:pr3.5}) again, we
find a bound of the form
\begin{eqnarray}
\label{eq:pr:22}
 \Bigg| \sum_{0<|k|< M}  \frac{C_{3}}{1+k^{2}}
 \frac{1}{r(|\zeta|^{2})} \sum_{\substack{ 
 \zeta' \in \Z^2 \\ |\zeta'| = |\zeta|}} 
 e^{i k \theta} \Bigg| \le 
  \max_{0<|k|<M} \;\; 
   \left|  \frac{ S(|\zeta|^{2},k)}{r(|\zeta|^{2})} 
     \right|
      \cdot
  \sum_{0<|k|<M} \frac{C_{3}}{1+k^{2}} \nonumber \\
            & &
\end{eqnarray}

Adding the error
terms~(\ref{eq:pr1.5}),~(\ref{eq:pr2}),~(\ref{eq:pr:21})
and~(\ref{eq:pr:22})  gives 
\begin{eqnarray}
\nonumber  \lefteqn{  | Q(f,f)(v) - Q^{h}(f^{h},f^{h})(v)|}&&
\hspace*{0.8\textwidth}\\ &\le& 
      \frac{C_{1}}{R^{2}} + C_{2} R^{2} h  +
      \frac{R^{2}  C_{4}}{M} + 
           C_{3} (2h)^{2}\max_{0<|k|<M} \sum_{\zeta\in Z_{h,R}} 
            \big|\frac{S(|\zeta|^{2},k)}{r(|\zeta|^{2})}\big|\nonumber
      \\
\label{eq:25_n} \
\end{eqnarray}
In the sum on the right hand side,
\begin{eqnarray*}
 \sum_{\zeta\in Z_{h,R}} 
      \left|  \frac{ S(|\zeta|^{2},k)}{r(|\zeta|^{2})} 
     \right|       =  \sum_{n<(R/h)^{2}}
            \big|S(n,k)\big| \,,
\end{eqnarray*}
and this can be estimated by using Proposition~\ref{prop7} with $X=(R/h)^{2}$.
To do this, we must require that 
\begin{equation}
  \label{eq:cc}
   R/h > \exp\left(\log(M)^{2}/b
   \right) 
\end{equation}
 for some positive constant $b$. Then there is a constant
$C_{5}$ such that 
\begin{eqnarray*}
 \sum_{n<(R/h)^{2}}\big|S(n,k)\big| \le C_{5}(\left(\frac{R}{h}\right)^{2}
 \exp\left(-\big(1-\frac{2}{\pi}) 
 \frac{\log\left( (R/h)^{2}\right)}{\left(\log M\right)^{2}}   \right) \,.
\end{eqnarray*}

The last term in~(\ref{eq:25_n}) will always be the dominating one,
and at this point, it does not give much to try to optimise the
choices of $R$, $M$ and $h$. Hence  to achieve an error of magnitude
$\varepsilon$ we
\begin{enumerate}
\item take $R = \sqrt{4 C_{1} / \varepsilon}$,
\item observe that we must have $h < \varepsilon / (4 R^{2} C_{2}) =
  \varepsilon^{2} / (4 C_{1}  C_{2})$,
\item choose $M = 4 R^{2} C_{4} / \varepsilon = 
           64  C_{1} C_{4} / \varepsilon^{2} $.
\end{enumerate}
With these choises of $R$ and $M$, the last term can then be bounded
by
\begin{eqnarray}
  4 C_{3} C_{5} \frac{4 C_{1}}{\varepsilon } 
  \exp\left(-\big(1-\frac{2}{\pi}) \log
 \frac{\log(4 C_{1} /(\varepsilon h^{2}) )}{\left(\log (64  C_{1}
  C_{4} / \varepsilon^{2})\right)^{2} }   \right)\,,
\end{eqnarray}
which converges to zero when $h\rightarrow 0$, and so there is an $h$
so small that also the last term in~(\ref{eq:25_n}) is smaller than
$\varepsilon/4$. We see that in order to achieve an error of
maginitude $\varepsilon$, one must take $h$  very small: $h = o\left( \exp(-
2\left(\log \varepsilon
\right)^{2}\varepsilon^{-2/(1-\frac{2}{\pi})})\right)$ (note
that~(\ref{eq:cc}) is then satisfied).

\end{proof}

\section{Some examples and remarks}
\label{sec:remarks}

From a numerical point of view, the discretization discussed
above would be far too costly:
a discrete velocity model with $N$
velocities would at least correspond to a computational cost of
$O(N)$ per time step, because one needs to compute a value
for each velocity. When the collision term is computed by the
sum~(\ref{eq:pr0}), the cost is  $ O(N^{2})$ times some
  logarithmic factor of $N$ (which comes from the summation over the
  points on the circles). And the calculation above showed that $N$
  grows exponentially in terms of the accuracy, $N \sim \frac{1}{h} >>
  \exp(\varepsilon^{-c}) $ for some positive constant $c$.

However, rather than estimating the computational cost in terms of the
number of discretization points used, it is more relevant to give the
cost in terms of the desired accuracy, given that the discretization
points are used in an optimal way. The discussion
around~(\ref{eq:pr0}) suggests that one can 
reduce the computational cost considerably without compromising  the
order of accuracy. The poor rate of convergence is due to the
approximation of $G_{v}(w)$.
Generalizing the formula~(\ref{eq:pr0}) slightly, we can
write
\begin{eqnarray}
\label{eq6.1}
    \int_{\R^{2}} G_{v}(w) \,dw \sim  
        \frac{1}{\rho_{h}} \sum_{\zeta\in Z_{h}} G_{v}(h\zeta)
\end{eqnarray}
where $\rho_{h}$ is the local density of $Z_{h}$. For $Z_{h} =
\{\zeta\in\Z^{2}\;\; \mbox{s.t.} \;\; |h\zeta| \le R\}$, one has
$\rho_{h} = h^{-2}$. Of course, even more generally  one could take a
local density which is not constant. 

The procedure for constructing a DVM would then be
\begin{itemize}
\item Choose a density $\rho_{h}$ so that the sum~(\ref{eq6.1}) is
  approximated to the desired order.
\item Choose $h$ so small that there exist a set $Z_{h}$ with this
  density so that for all $\zeta\in Z_{h}$  the angular integral is
  well approximated by the sum.
\end{itemize}

For such a model, the computational cost for each velocity would be
of the order $\varepsilon^{-3}$ (this estimate is based on the
assumption that the cost of evalutating the angular integral is
$\varepsilon^{-1}$, and that the number of velocities is
$O(\varepsilon^{-2})$; lower cost can be acheived if higher order formulas 
are used for approximating the integrals). The problem remains, that a
very large number of velocities are needed, and hence the total
computational cost is still excessive. A more challenging task would
be to dilute not only the set $Z_{h}$, but to  choose in a
systematic way subsets  $U_{h}\subset h\Z^{2}$ for the discrete velocity model,
so that $Q(f,f)(v)$ would be well approximated for all $v\in U_{h}$,
and to do this in a way that does not require too large tables for
storing all possible collisions.

In the last part of this paper, we wish to illustrate the distribution
of good radii. We then consider  $\zeta =
(\zeta_{1},\zeta_{2}) \in \Z^{2} \;\;
\mbox{s.t.}\;\; 0 \le \zeta_{i} \;\; (i=1,2)\;\; |\zeta|<20000
\}$. This is an 
extremely large set of points, which corresponds to a huge 
number  of velocities (the $ O(N^{2})$ factor would in
this case be of the order $10^{17}$, which is of course absurd)
 
\begin{figure}[htbp]
  \centering
  \includegraphics[width=12cm]{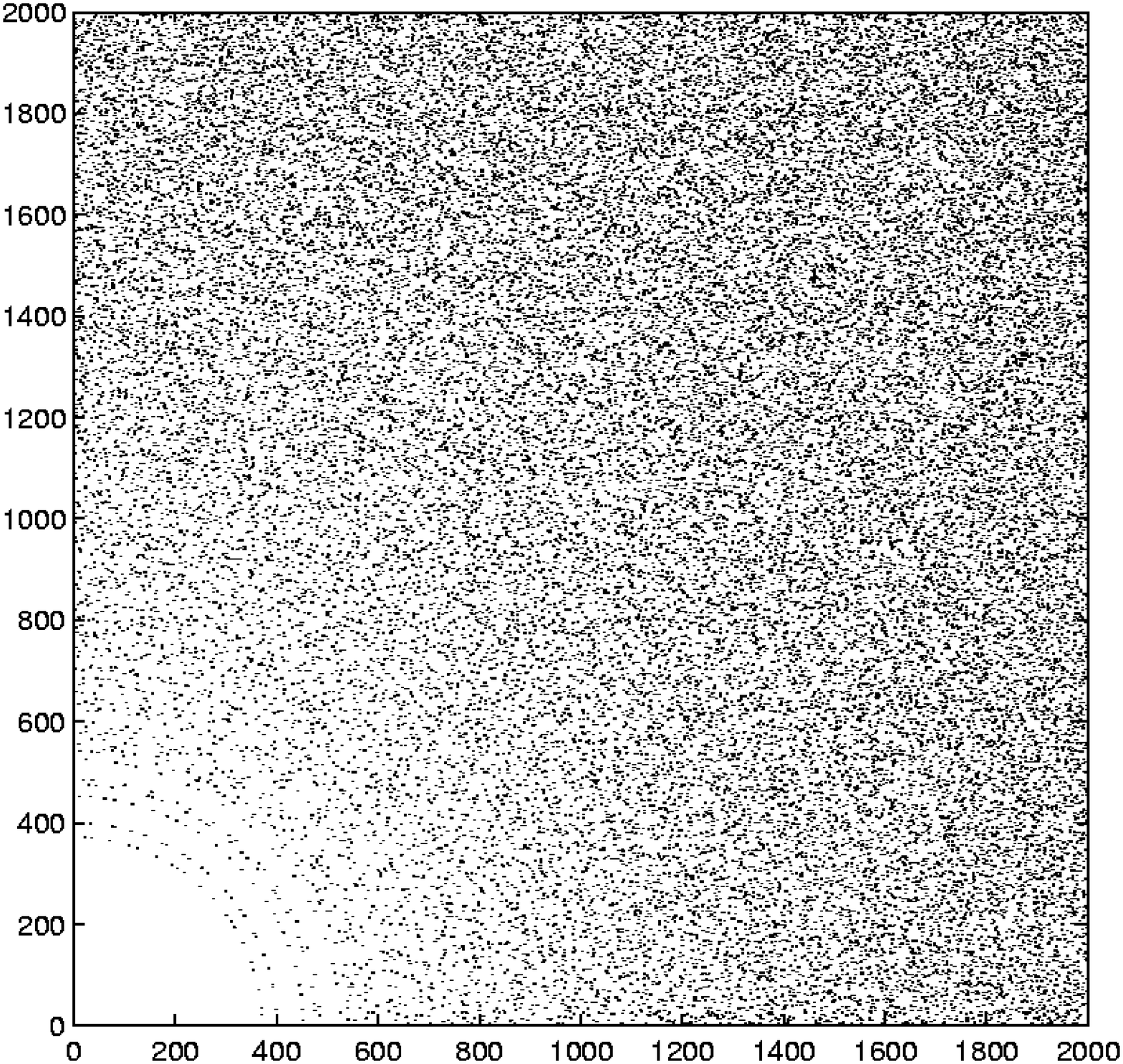}
  \caption{Lattice points such that circles containing these points,
  contain at least 72 lattice points  }
  \label{fig_1}
\end{figure}

Among the circles with radii $|\zeta|$ in this set, the largest number
of points on one circle, is 384.
In Fig.~\ref{fig_1}, we show all points $\zeta =
(\zeta_{1},\zeta_{2})$ with $0<\zeta_{i}<2000$, such that the circle
passing through $\zeta$  has more than 72 points. There are 36163
points in this set.  This is a small fraction of the total number of
integer points, but   they are seemingly well distributed, except near
the origin.

\begin{figure}[htbp]
  \centering
  \includegraphics[width=12cm]{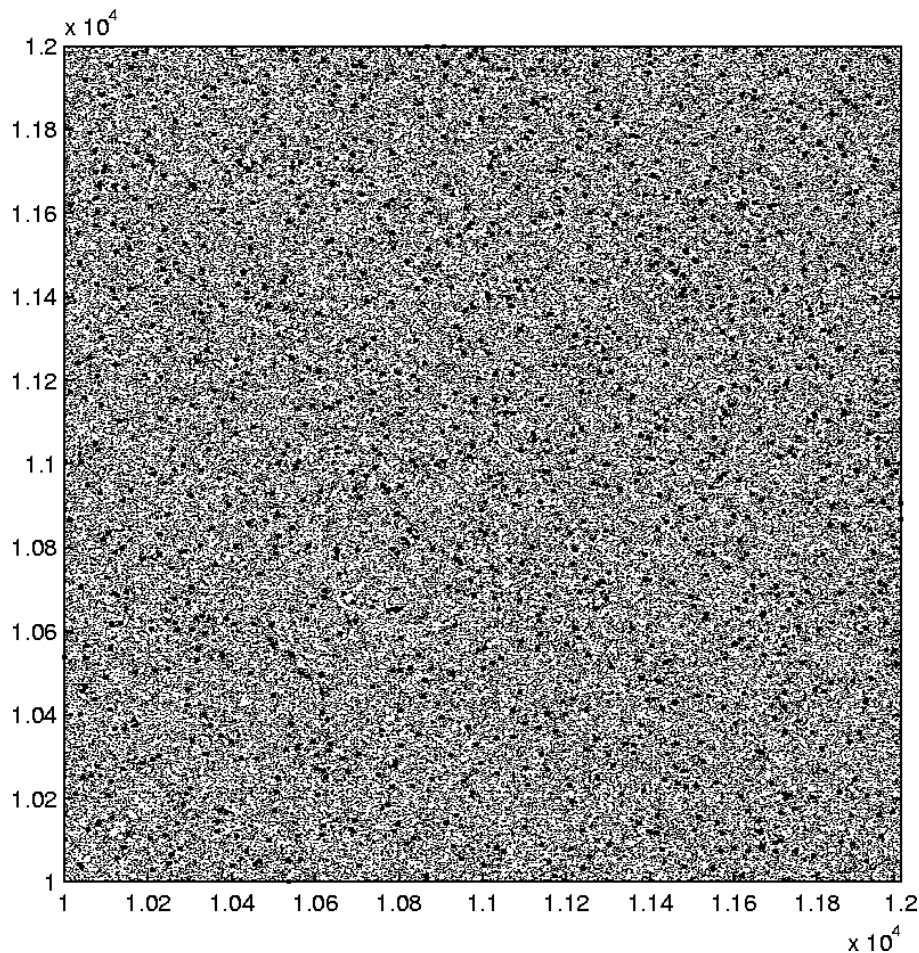}
  \caption{Lattice points such that circles containing these points,
  contain at least 72 lattice points (small dots), or at least 192
  points (the larger dots) }
  \label{fig_2}
\end{figure}

Figure~\ref{fig_2} shows points in the range $10000\le \zeta_{i} \le
12000$. Here the small dots denote points on circles having at least
72 points, and the larger dots denote points on circles with at least
192 points (there are 141562 and 1120 points respectively in these
sets).

\bigskip

\noindent{\bf Acknowledgment:} We would like to thank A. Bobylev,
J. Brzezinski, A. Heintz, and Z. Rudnick for useful discussions.

\end{document}